\documentstyle[mssymb]{article}
\def\MW{Mordell-\mbox{\kern-.1em}Weil}
\newcommand{\be}{\begin{equation}}
\newcommand{\ee}{\end{equation}}
\newcommand{\bea}{\begin{eqnarray}}
\newcommand{\eea}{\end{eqnarray}}
\newcommand{\ra}{\rightarrow}
\newcommand{\naive}{na\"\i ve}
\newcommand{\Ht}{{\hat h}}

\newcommand{\Erho}{{E_8^\rho}}

\newcommand{\eps}{\epsilon}
\newcommand{\veps}{\varepsilon}
\newcommand{\om}{\omega}
\newcommand{\vrho}{\varrho}
\def\bmu{{\mbox{\boldmath$\mu$}}}
\font\atebf=cmbx8
\newcommand{\sF}{{\hbox{\atebf F}}} 
\newcommand{\sQ}{{\hbox{\atebf Q}}} 
\font\svnbf=cmbx7
\newcommand{\subZ}{{\hbox{\atebf Z}}} 
\newcommand{\Aff}{{\bf A}}
\newcommand{\A}{{\cal A}}
\newcommand{\E}{{\cal E}}
\newcommand{\C}{{\bf C}}
\newcommand{\sfC}{{\sf C}}
\newcommand{\F}{{\bf F}}
\newcommand{\M}{{\cal M}}
\newcommand{\Q}{{\bf Q}}
\newcommand{\R}{{\bf R}}
\newcommand{\SS}{{\cal S}}
\newcommand{\Z}{{\bf Z}}
\renewcommand{\P}{{\bf P}}
\newcommand{\Six}{{({\sP^1\atop 6_{\phantom1}})}}
\newcommand{\SIX}{{\bf P}^1_{\!(6)}}
\newcommand{\Five}{{({\sA^1\atop 5_{\phantom1}})}}

\font\svnbf=cmbx7
\newcommand{\sP}{{\hbox{\svnbf P}}} 
\newcommand{\sA}{{\hbox{\svnbf A}}} 
\newcommand{\sZ}{{\hbox{\svnbf Z}}} 

\newcommand{\NS}{{\rm NS}}

\newcommand{\Aut}{\mathop{\rm Aut}}

\newcommand{\Gal}{\mathop{\rm Gal}}

\newcommand{\Hom}{\mathop{\rm Hom}}
\newcommand{\GL}{\mathop{\rm GL}}
\newcommand{\PGL}{\mathop{\rm PGL}}
\newcommand{\SL}{\mathop{\rm SL}}

\newcommand{\Sp}{\mathop{\rm Sp}}
\newcommand{\GSp}{\mathop{\rm GSp}}
\newcommand{\PSp}{\mathop{\rm PSp}}
\newcommand{\PGSp}{\mathop{\rm PGSp}}
\newcommand{\SO}{\mathop{\rm SO}}
\newcommand{\U}{\mathop{\rm U}}
\newcommand{\PU}{\mathop{\rm PU}}
\newcommand{\SU}{\mathop{\rm SU}}
\newcommand{\SigU}{\mathop{\Sigma\rm U}}

\newcommand{\0}{^{\phantom0}}
\renewcommand{\baselinestretch}{1.3}
\begin{document}

\begin{center}
\begin{large}
The identification of three moduli spaces\\
\vspace*{2ex}
Noam D. Elkies\\
\end{large}
May, 1999

\end{center}

\vspace*{4ex}

{\small
{\bf Abstract.}  It is one of the wonderful ``coincidences''
of the theory of finite groups that the simple group~$G$\/
of order 25920 arises as both a symplectic group in characteristic~3
and a unitary group in characteristic~2.  These two realizations
of~$G$\/ yield two $G$-covers of the moduli space $\Six$ of
configurations of six points on the projective line modulo $\PGL_2$,
via the 3- and 2-torsion of the Jacobians of the double and triple
cyclic covers of $\P^1$ branched at those six points.  Remarkably
these two covers are isomorphic.  This was proved over~$\C$
by transcendental methods in~\cite{Janus}.  We give an algebraic
proof valid over any field not of characteristic~2 or~3 that
contains the cube roots of unity.  We then explore the connection
between this $G$-cover $\SS$\/ of~$\Six$ and the elliptic surface
$y^2 = x^3 + {\rm sextic}(t)$, whose \MW\ lattice is $E_8$ with
automorphisms by a central extension of~$G$. 
}

\vspace*{2ex}

{\bf 0. Introduction.}  The moduli spaces of the title all
cover the moduli space $\Six$ of unordered sextuples
of distinct points on $\P^1$ modulo the action of $\PGL_2$,
or equivalently of sextic polynomials $S(t)$ without repeated roots
modulo the action of $\GL_2$.  [A quintic counts as a sextic with
a root at infinity; polynomials of lower degree would have repeated
roots at infinity and are thus excluded.]  To a configuration of six
distinct points on $\P^1$ are associated two curves $C: u^2 = S(t)$
and $C': v^3 = S(t)$ of genus~2 and~4, which are cyclic covers
of $P^1$ of degrees~2 and~3 branched at those six points.
Two of our covers of~$\Six$ are obtained by adding full level-3
structure to~$C$\/ and full level-2 structure to~$C'$; that is,
they are the moduli of genus-2 curves $C$\/ with a choice of
generators for the 3-torsion subgroup of the Jacobian $J(C)$,
and of curves $C'$ with a choice of generators for the 2-torsion
subgroup of $J(C')$.  The Weil pairing gives these torsion subgroups
the structure of four-dimensional spaces over the finite fields
$\F_3$ and $\F_4$, with respectively a symplectic and a unitary
structure; thus the corresponding moduli spaces are normal covers
of $\Six$ with Galois groups $\PSp_4(\F_3)$ and $\SU_4(\F_2)$.
It is one of the wonderful ``coincidences'' of the theory of
simple finite groups that these two groups of order $25920$ coincide;
let $G$\/ be the finite group isomorphic with both $\PSp_4(\F_3)$
and $\SU_4(\F_2)$.  (See page 26 of the ATLAS \cite{ATLAS} for this
identification and for and further properties of~$G$\/ that we
shall use.)  It is a remarkable fact that the two $G$-covers
of $\Six$ are also isomorphic.  This is proved in~\cite{Janus}
(see also \cite[Ch.5]{Hunt}) by transcendental methods:
first replace $\Six$ by the moduli space $\SIX$
of ordered sextuples of distinct points, which is an $S_6$ cover
of~$\Six$; then regard both $G$-covers of $\SIX$ over~$\C$
as quotients of the complex 3-ball by arithmetic groups acting freely;
and prove that they are isomorphic by computing enough invariants.
Our two $G$\/-covers of $\Six$ are then quotients of these $G$\/-covers
of $\SIX$ by the same action of $S_6$, are thus isomorphic as well.
Of course this leaves completely mysterious the algebraic meaning
of the identification between the two moduli spaces.

We obtain this identification algebraically via a third moduli
space $\SS$: the space of sextic polynomials $S(t)$ together with
all representations of $S$\/ as the difference between the square
and the cube of polynomials of degree at most~3 and~2, which we
shall call {\em minimal}\/ representations of~$S(t)$ as $y^2-x^3$.
By counting parameters we may surmise that the forgetful
map from $\SS$\/ to $\Six$ is a Galois cover,
but it is not at all clear what the group should be.
We prove that the group is~$G$\/ by constructing algebraic maps
from $\SS$\/ to the two geometric $G$\/-covers of~$\Six$
and showing that the maps are isomorphisms.  This of course yields
an algebraic isomorphism between the $J(C)[3]$ and $J(C')[2]$
moduli spaces.  We then give a further interpretation of~$G$\/
as the group of automorphisms of a $\Z[e^{2\pi i/3}]$ lattice
in $\C^4$ isomorphic as a Euclidean lattice with~$E_8$,
which arises here as the \MW\ lattice of the rational
elliptic surface
\be
\E = \E_S: y^2=x^3+S(t).
\label{Surface}
\ee
The 240 minimal vectors of this lattice correspond $2:1$
with the 120 odd elements of order~$2$ in the Jacobian of $C'$;
they correspond $3:1$ with the 80 nontrivial elements of order~$3$
in the Jacobian of~$C$\/; and they correspond $6:1$ with the 40
representations of $S$\/ as the difference between the cube of
a quadratic and the square of a cubic polynomial.  Even the
fact that there are always 40 such representations is far
from well known, though it turns out that Clebsch \cite{Clebsch}
had already obtained this enumeration and, in collaboration with
Jordan, also its relation with what we now call $J(C)[3]$.

{\em
After this paper was largely completed I found that the identification
of the $J(C)[3]$ and $J(C')[2]$ moduli spaces via the minimal solutions
of $y^2=x^3+S(t)$ was already obtained by van~Geemen in~\cite{Geemen},
and in essentially the same way.  He was led to it via the study of
theta functions.  The interpretation of this result in terms of the
\MW\ lattice of~(\ref{Surface}), and the connections with other
geometric ideas given or announced here, still appear to be new.
}

{\bf Acknowledgements.}  Thanks to Daniel Allcock for bringing
\cite{Hunt} to my attention, and to Allcock and Joe Harris for
helpful discussions.  Jordan Ellenberg read an earlier version
of this paper and found several instances of confusing or mistyped
prose whose correction markedly improved the exposition.

This work was made possible in part by funding from the
David and Lucile Packard Foundation.

{\bf 1. The $\M_2\0$ and $\M_4^\vrho$ pictures of $\Six$
and their $G$-covers.}
Let $k$\/ be a field not of characteristic~$2$ or~$3$
containing a cube root of unity~$\rho$.
To six distinct points $t_1,\ldots,t_6$ on $\P^1(\bar k)$
permuted by $\Gal(\bar k/k)$ we associate two curves defined
over~$k$\/: the genus~2 curve
\be
C: u^2 = S(t) = \eps\prod_{j=1}^6 (t-t_j)
\label{C}
\ee
up to quadratic twist $\eps\in k^*/{k^*}^2$, and the genus-4
superelliptic curve (called a ``Picard curve'' in~\cite{Janus})
\be
C': v^3 = S(t) = \veps\prod_{j=1}^6 (t-t_j)
\label{C'}
\ee
up to cubic twist $\veps\in k^*/{k^*}^3$.
If some $t_j=\infty$ the corresponding factor $t-t_j$
is replaced by~$1$.  Except for the twists $\eps,\veps$,
these curves do not depend on the choice of coordinate on $\P^1(k)$:
changing $t$\/ to $(at+b)/(ct+d)$ changes $\prod_j (t-t_j)$ to
\be
\prod_{j=1}^6 \left( \frac{at+b}{ct+d} - \frac{at_j+b}{ct_j+d} \right)
= \left(\frac{ad-bc}{ct+d}\right)^{\!6}
  \prod_{j=1}^6 \frac{t-t_j}{ct_j+d},
\label{t_change}
\ee
in which the factor $\bigl((ad-bc)/(ct+d)\bigr)^6$ may be absorbed
into $u^2$ or $v^3$, while the constant factor $\prod_{j=1}^6 (t_j+d)$
is absorbed into $\eps$ or $\veps$.  The reader may check that in the
special cases $t_j=\infty$ or $ct_j+d=0$ the formula (\ref{t_change})
still holds with our interpretation of a factor ``$t-\infty$'' as~$1$,
and the factor $ct_j+d$ replaced by $-c$ if $t_j=\infty$ and by
$-(at_j+b)$ if $ct_j+d=0$.

It is well known that every curve of genus~$2$ is of the form $C$\/
for some choice of $t_j$ and $\epsilon$, uniquely determined by
the curve up to the action of $\PGL_2(k)$, and that every principally
polarized abelian surface which is not a product of elliptic curves is
the Jacobian of such a curve.  In other words, it is known that $\Six$
is identified with the moduli space $\M_2$ of curves of genus~$2$,
and with the open subset of the moduli space $\A_2$ of principally
polarized abelian surfaces, namely the subset parametrizing
indecomposable surfaces.

It is also known, though not as widely, that the curves $C'$
yield another pair of moduli interpretations of $\Six$.
Just as the genus-2 curves $C$\/ have the hyperelliptic involution
$(t,u)\leftrightarrow(t,-u)$, the genus-4 curves have an automorphism
of order~$3$ defined by
\be
\vrho(t,v) = (t,\rho v).
\label{rho}
\ee
The fixed points of this involution are the six points above $t=t_j$.
By the Riemann-Hurwitz formula, an order-3 automorphism of any curve
of genus~4 has six fixed points if and only if the quotient curve has
genus~0.  Now there are two kinds of genus-4 curves cyclically covering
$\P^1$ with degree~$3$: the curves $C'$, and curves of the form
\be
v^3 = \prod_{i=1}^3 \frac{t-t_i}{t-t_j}.
\label{notC'}
\ee
The two kinds of curves are distinguished by the action of $\vrho$ on
the four-dimensional space of holomorphic differentials on the curve.
In both cases the fixed subspace is necessarily trivial because
it is the space of holomorphic differentials on the quotient curve.
But the action of $\vrho$ on $H^1(C')$ is diagonalized by
the basis $(dt/v^2, t\,dt/v^2, t^2\,dt/v^2, dt/v)$, the first three
of whose vectors have eigenvalue $\rho$ and the last has eigenvalue
$\bar\rho$; whereas the $\rho$- and $\bar\rho$-eigenspaces of the
holomorphic differentials on (\ref{notC'}) both have dimension~$2$.
It follows that over an algebraically closed field the curves $C'$
are precisely those genus-4 curves with an order-3 automorphism $\vrho$
whose $1$-, $\rho$-, and $\bar\rho$-eigenspaces have dimensions $0,3,1$.
When $k$\/ is not algebraically closed we must also impose the condition
that the quotient curve have a rational point (and thus be identified
with $\P^1$).\footnote{
  If we worked over a field~$k$\/ not containing $\rho$ we would
  also have to require that the involution in $\Gal(k(\rho)/k)$ take
  $\vrho$ to $\vrho^{-1}$.
  }
We thus have a straightforward identification of $\Six$
with a space we shall call $\M_4^\vrho$, the moduli space of
genus-4 curves $C'$ with an order-3 automorphism $\vrho$
whose action on $H^1(C)$ has eigenvalues $\rho,\rho,\rho,\bar\rho$.
Less obvious, but still true, is that the Jacobian map $\M_4\ra\A_4$
identifies $\M_4^\vrho$ with an open subset of the moduli space
$\A_4^\vrho$ of what we shall call ``$\vrho$-fourfolds'':
principally polarized abelian fourfolds with an
order-3 automorphism $\vrho$ acting on the tangent space
with eigenvalues $\rho,\rho,\rho,\bar\rho$.  This is seen by computing
that $\A_4^\vrho$ has dimension~$3$, the same as the dimension of
$\M_4^\vrho$.~\cite{Janus}

Now each of $\A_2\0$ and $\A_4^\vrho$ has natural arithmetic covers:
the moduli spaces of principally polarized abelian surfaces or
$\vrho$-fourfolds with additional (torsion) structure.  We will
thus obtain two families of natural covers of $\Six$.  For instance,
for each prime $p$ the moduli space of principally polarized abelian
surfaces with full level-$p$ structure --- that is, with a choice
of generators for their $p$-torsion --- is a Galois cover of $\A_2$
with Galois group $\PSp_4(\F_p) = \Sp_4(\F_p)/\{\pm1\}$
if the ground field contains the $p$-th roots of unity.
[Not all of $\GL_4(\F_p)$ because the Galois group must respect
not only the group structure but also the Weil pairing on the
$p$-torsion; we divide by $\{\pm1\}$ to account for
the quadratic twists.]  When $p=2$, a full level-2 structure
on $J(C)$ is just an ordering of the six Weierstrass points,
because the 15 nontrivial 2-torsion points are represented by
the differences between pairs of Weierstrass points; thus we
obtain the cover of $\Six$ by the moduli space $\SIX$ of
{\em ordered}\/ sextuples of distinct points on $\P^1$ modulo $\PGL_2$.
Of course the Galois group of this cover is the symmetric group $S_6$,
so we have recovered the identification of this group with the
symplectic group $\Sp_4(\F_2)$.  Likewise we obtain Galois covers
of $\A_4^\vrho$ from the $p$-torsion points of $\vrho$-fourfolds,
with the Galois group this time depending on whether $p$ is ramified,
inert or split in $\Q(\rho)$.  For the
ramified prime~$3$ we obtain an intermediate cover from the 3-torsion
points in the kernel of $\sqrt{-3}=\vrho-\bar\vrho$.  These constitute
a 4-dimensional space over $\F_3$; the Weil pairing yields a quadratic
from on this space, taking a $(\sqrt{-3})$-torsion point~$P$\/ to
the pairing of $P$\/ with any of the 3-torsion points $P'$ such that
$P=\sqrt{-3} P'$.  This quadratic form turns out to have Arf
invariant~$1$, so we obtain a cover of $\A_4^\vrho$ with Galois group
$\SO^-_4(\F_3)$.  Again the $(\sqrt{-3})$-torsion points on $J(C')$
are generated by the differences between the six points with $t=t_j$
(note that these are fixed by $\vrho$, so in the kernel of
$1-\rho=\bar\rho\sqrt{-3}$), so once more we find the cover of
$\Six$ by $\SIX$, whose Galois group $S_6$ is this time
identified with $\SO^-_4(\F_3)$.  (See \cite[p.4]{ATLAS}
for the realizations of $S_6$ as linear groups in characteristics
2 and~3.)

The next cases are $p=3$ for $\A_2$ and $p=2$ for $\A_4^\vrho$.
We already know that the former yields a cover of $\A_2$ with
Galois group $\PSp_4(\F_3)=G$.  As to the latter, the even prime
is inert in $\Q(\rho)$, so the 2-torsion points of $J(C')$ have
the structure of a four-dimensional vector space over~$\F_4$.
The Weil pairing is an alternating $\F_2$-bilinear map on that
space consistent with the $\F_4$ structure (that is, such that
$\langle P,P' \rangle = \langle \vrho P, \vrho P' \rangle$); this
yields a natural unitary form on $J(C')[2]$, namely
$\upsilon: P \mapsto \langle P,\vrho P\rangle$.  Thus the
Galois group for the cover of $\A_4^\vrho$ by the
moduli space of $\vrho$-fourfolds with a full level-2
structure is $\U_4(\F_2)/\F_4^* = \SU_4(\F_2)$ is again~$G$.
(We divide by the center $\F_4^*$ to account for the cubic twists;
since the dimension~4 is coprime to $3=\#\F_4^*$, the resulting group
$\PU_4(\F_2)$ is isomorphic with $\SU_4(\F_2)$.
So we have again obtained, from two different curves
$C,C'$ associated with a six-point configuration in $\P^1$ and
full level structures of different levels, covers of $\Six$
with the same Galois group.  But, unlike the case of the $S_6$-cover
of $\Six$ by $\SIX$, it is not immediately clear here that the
two $G$-covers of $\Six$ are the same.  This is what Hunt and Weintraub
proved transcendentally in~\cite{Janus}, and what we next prove
algebraically by identifying both spaces with the space of sextics
$S(t)$ together with all minimal solutions $(x,y)\in k[t]\times k[t]$
of $y^2=x^3+S(t)$.

\vspace*{4ex}

{\bf 2. Torsion divisors and $\SS$.}
The action of~$\GL_2$ on sextic polynomials $S(t)$ also respects
the minimal solutions of $y^2=x^3+S(t)$.
Indeed, substituting $(at+b)/(ct+d)$ for~$t$ in
$y^2=x^3+S(t)$ and multiplying by $(ct+d)^6$ yields
$y_1^2 = x_1^3 + S_1(t)$, where $x_1,y_1,S_1$ are the images
of~$x,y,S$\/ under $({{a\;b}\atop{c\;d}})\in\GL_2$.
Thus the set of minimal solutions of $y^2=x^3+S$\/
makes sense even if $S$\/ is not a sextic
polynomial but a $\GL_2$-orbit of such polynomials, or a point
on~$\Six$.  In the sequel we assume for simplicity that no $t_j=\infty$,
i.e.\ that $S(T)$ is of degree~6, not~5.  The argument can be readily
adapted to handle the case of a quintic; alternatively
we may make all the $t_j$ finite by first applying a $\PGL_2(k)$
transformation, which as already observed does not change the
problem.  Note that since $k\ni\rho$ and $k$\/ is not of characteristic
$2$ or~$3$ there are always more than six points in $\P^1(k)$ so
such a transformation must exist.

{\bf Theorem 1.}
{\sl
Let $y\in k[t]$ be a polynomial of degree at most\/~$3$ such that
$y^2-S(t)$ is a $k^*$-multiple of a cube in $k[t]$.  Then the rational
function $y(t)-u$ on~$C$ has divisor $3D$\/ for some nonprincipal
divisor~$D$.  Each nonzero 3-torsion element of $J(C)$ arises as
the class of such a divisor~$D$\/ for exactly one choice of~$y$.
}

{\sl Remark}\/: This theorem, and its following proof, turn out to
have already been given (in older terminology, and only over~$\C$)
in~\cite{Clebsch}.

{\sl Proof}\/: On $C$\/ we have $S(t)=u^2$, so $\alpha x^3=y^2-S(t)$
factors as $\alpha x^3=(y-u)(y+u)$.  Now $y,S$\/ have no common roots:
at such a root $x^3$ would vanish, so it would be a root of~$x$
as well; but then $S=y^2-\alpha x^3$ would vanish there to
order at least~$2$, which is not allowed.  Thus $y,u$ have
no common zeros on~$C$, and $u$ has a triple pole at each of
the points at infinity while $y$ has a pole of order at most~$3$
there, and neither $y$ nor $u$ has any pole at a point of~$C$\/
with $t\neq\infty$.  Thus the $y-u$ and $y+u$ have no common zeros
nor any finite poles, and at each point at infinity at least one of
$y\pm u$ has a triple pole while the other has a pole of order
$\leq 3$.  Now $x$ is regular except for poles of order at
most~$2$ at the points at infinity.  Since $(y+u)(y-u)=\alpha x^3$,
we conclude that the divisor of $y+u$ has the form $3D$\/ for some
divisor~$D$\/ on~$C$.  Since $3D$\/ is principal, the class
$[D]\in J(C)$ of~$D$\/ is a 3-torsion point on the Jacobian.
But $D$\/ cannot itself be principal.  Indeed, if it were,
$(y+u)^{1/3}$ would be a rational function on~$C$\/ of degree
at most~$2$.  But all such functions are in $k(t)$.
Thus $y+u$ would also in $k(t)$, which is absurd because
$y\in k(t)$ but $u\notin k(t)$.  Thus $[D]$ is a nonzero
3-torsion point in the Jacobian, as claimed.

Now let the divisor $K$\/ on~$C$\/ be the sum of the two points at
infinity; this is the divisor of the differential $dt/u$ on~$C$,
and is thus canonical.  Since $D$\/ includes at worst simple poles
at the points at infinity, $D+K$\/ is an effective, non-canonical
divisor of degree~$2$.  Conversely, let $T$\/ be a nonzero 3-torsion
point in $J(C)$.  Then $K+T$ is a non-canonical divisor class of
degree~$2$, and so by Riemann-Roch has a unique effective
representative, which we may call $K+D$.  Since $D\sim T$,
the divisor $3D$\/ is principal, say the divisor of $f\in(k(C))^*$.
Since $3K+3D$\/ is effective, this function $f$\/ has at most triple
poles at the points at infinity, and no other poles; thus it is of the
form $y_1(t)+cu$ where $y_1$ is a polynomial of degree at most~$3$
and $c\in k$.  If $c$ vanished then $f$\/ would be a
rational function of~$t$\/ which, when considered as a function
on~$C$, would have all its zeros and poles of order divisible by~$3$.
But then the same would be true of~$f$\/ as a function on~$\P^1$,
because the cover $t:C\ra\P^1$ has degree prime to~$3$.  Thus $f$\/
would be a $k^*$-multiple of the cube of a rational function
--- and this rational function would have divisor~$D$, which would
therefore be principal, contradicting the hypothesis that $T\neq0$.
Thus $f=y_1+cu$ for some $c\in k^*$, and multiplying by
$c^{-1}$ we find the unique function of the form $y+u$ whose
divisor is~$3D$\/ for some $D\sim T$.  It remains to prove that
$y^2-S(t)$ is a $k^*$-multiple of a cube.  Since the divisor
of $y+u$ is divisible by~3, the same is true of the divisor
of $y-u$, which is the image of~$y+u$ under the hyperelliptic
involution.  But then the same is true of $(y-u)(y+u)=y^2-S(t)$,
considered as a function on~$C$.  As noted already this means
that the divisor of $y^2-S(t)$ considered as a function of~$t$\/
is also divisible by~$3$, so $y^2-S(t)$ is indeed of the form
$\alpha x^3$ for some $x\in k[t]$ and we are done.

{\sl Corollary}: i) \cite{Clebsch}
There are $3^4-1=80$ polynomials $y\in\bar k[t]$ of degree
at most~$3$ such that $y^2-S(t)$ is a cube in~$\bar k[t]$.\\
ii) $\SS$ is a $G$-cover of~$\Six$, isomorphic with the moduli
space of curves of genus~$2$ with full level-3 structure.~~$\Box$

Much the same analysis applies to the curves $C'$.  The one difference
here is that the divisor $K_0$, consisting of the sum of the points
(now three of them) at infinity, is no longer canonical: the
differential $dt/v^2$ has divisor $2K_0$, not $K_0$.  Thus $K_0$
is a distinguished {\em semicanonical}\/ divisor, a.k.a.\ theta
characteristic.  (As with the divisor of~$K$\/ on~$C$, the divisor
$K_0$ depends on the choice of coordinate~$t$\/ on $\P^1$,
but the linear equivalence class of $K_0$ does not: if we
used instead a coordinate whose pole is at $\tau\neq\infty$
our divisor would be the fiber of~$\tau$ under the map $t:C'\ra\P^1$,
and the difference between these two $K_0$'s would be the divisor
of the rational function $t-\tau$.)  Using $K_0$ we may identify
2-torsion elements~$T$\/ of the Jacobian with theta characteristics
$K_0+T$.  Recall that there is an affine-linear quadratic form
on the theta characteristics which gives the parity of the space
of sections; for genus~4, there are $(2^8+2^4)/2=136$ even and
$(2^8-2^4)/2=120$ odd theta characteristics.  The divisor $K_0$
has the two-dimensional space of sections generated by~$1$ and~$t$,
so it is an even theta characteristic. 
None of the other 135 even thetas~$\Theta$ has any nonzero sections:
if it did then by parity it would have at least~$2$ independent ones,
whose ratio would be a degree-3 map to~$\P^1$ not of the form
$(at+b)/(ct+d)$; but there are no such rational functions on~$C'$.
The same argument shows that each odd theta has only a one-dimensional
space of sections, and thus has a unique effective representative.

We are now ready to state the analogue of Theorem~1 for $C'$:

{\bf Theorem 2.}
{\sl
Let $x\in k[t]$ be a polynomial of degree at most\/~$2$ such that
$x^3-S(t)$ is a $k^*$-multiple of a square in $k[t]$.  Then the rational
function $x(t)+v$ on~$C'$ has divisor $2D$\/ for some nonprincipal
divisor~$D$\/ such that $[D+K_0]$ is an odd theta characteristic.
Each odd theta characteristic arises in thus way for exactly one
choice of~$x$.
}

{\sl Proof}\/: On $C'$ we have $S(t)=v^3$, so $\beta y^2 = x^3+S(t)$
factors as $\beta y^2=(x+v)(x+\rho v)(x+\bar\rho v)$.  Now $x,S$\/
have no common roots: at such a root $y^2$ would vanish, so it would
be a root of~$y$ as well; but then $S=\beta y^2 - x^3$ would vanish
there to order at least~$2$, which is not allowed.  Thus $x,v$ have
no common zeros on~$C'$, and $v$ has double poles at each of the
points at infinity while $x$ has a pole of order at most~$2$
there, and neither $x$ nor $v$ has any pole at a point of~$C'$
with $t\neq\infty$.  Thus, of the three factors $x+v$, $x+\rho v$,
and $x+\bar\rho v$ of $\beta y^2$, none has a pole except at one of
the points at infinity; each has poles at infinity of order at
most~$2$; at each point at infinity at most one factor may have a pole
of order $<2$; and no two have a common zero.  Now $y$ is regular
except for poles of order at most~$3$ at the points at infinity.
Since $(x+v)(x+\rho v)(x+\bar\rho v) = \beta y^2$, we conclude that the
divisor of $x+v$ has the form $2D$\/ for some divisor~$D$\/ on~$C'$.
Since $2D$\/ is principal, the class $[D]\in J(C')$ of~$D$\/
is a 2-torsion point on the Jacobian.  But $D$\/ cannot itself
be principal.  Indeed, if it were, $(x+v)^{1/2}$ would be a rational
function on~$C'$ of degree at most~$3$.  But all such functions are
in $k(t)$.  Thus $x+v$ would also in $k(t)$, which is absurd because
$x\in k(t)$ but $v\notin k(t)$.  Thus $[D]$ is a nonzero 2-torsion
point in the Jacobian, as claimed.  Moreover $[D+K_0]$ is an odd
theta characteristic, because $D+K_0$ is an effective semicanonical
divisor not equivalent to~$K_0$.

Conversely, let $T$\/ be a nonzero 2-torsion point in $J(C')$
such that $K_0+T$\/ is an odd theta characteristic, and let
$K_0+D$\/ be the unique effective divisor in the class of $K_0+T$.
Since $D\sim T$, the divisor $2D$\/ is principal, say the divisor
of $f\in(k(C'))^*$.  Since $2K+2D$\/ is effective, this function $f$\/
has at most double poles at the points at infinity, and no other poles;
thus it is of the form $x_1(t)+cu$ where $x_1$ is a polynomial of
degree at most~$2$ and $c\in k$.  If $c$ vanished then $f$\/ would be a
rational function of~$t$\/ which, when considered as a function
on~$C'$, would have all its zeros and poles of order divisible by~$2$.
But then the same would be true of~$f$\/ as a function on~$\P^1$,
because the cover $t:C'\ra\P^1$ has degree prime to~$2$.  Thus $f$\/
would be a $k^*$-multiple of the square of a rational function
--- and this rational function would have divisor~$D$, which would
therefore be principal, contradicting the hypothesis that $T\neq0$.
Thus $f=x_1+cv$ for some $c\in k^*$, and multiplying by
$c^{-1}$ we find the unique function of the form $x+v$ whose
divisor is~$3D$\/ for some $D\sim T$.  It remains to prove that
$x^3+S(t)$ is a $k^*$-multiple of a square.  Since the divisor
of $x+v$ is divisible by~2, the same is true of the divisors of
$x+\rho v$ and $x+\bar\rho v$, which are the images of $x+v$
under $\vrho$ and $\vrho^2$.  But then the same is true of
$(x+v)(x+\rho v)(x+\bar\rho v) = x^3 + S(t)$ considered as
a function on~$C'$.  As noted already this means
that the divisor of $x^3+S(t)$ considered as a function of~$t$\/
is also divisible by~$2$, so $x^3+S(t)$ is indeed of the form
$\beta y^2$ for some $y\in k[t]$ and we are done.

{\sl Corollary}: i) There are $120$ polynomials $x\in\bar k[t]$ of
degree at most~$2$ such that $x^3+S(t)$ is a square in~$\bar k[t]$.\\
ii) $\SS$ is a $G$-cover of~$\Six$, isomorphic with the moduli
space of curves~$C'$ with full level-2 structure.~~$\Box$

Combining the Corollaries to Theorems 1 and~2 we obtain

{\bf Theorem 3.} {\sl
i) Any polynomial $S(t)$ of degree $5$ or~$6$ without repeated
roots over an algebraically closed field not of characteristic
$2$ or~$3$ can be written as $y^2-x^3$ for $240$ pairs of
polynomials $x(t),y(t)$ of degree at most\/~$2$ and\/~$3$
respectively.\\
ii) Over any field not of characteristic $2$ or~$3$ that contains
the cube roots of unity, the moduli space of curves $C$ with
full level-3 structure and the moduli space of curves $C'$
with full level-2 structure are isomorphic $G$-covers of~$\Six$.
}

{\sl Proof}\/: i) This follows from part (i) of the Corollary
to either Theorem~1 or Theorem~2 since $240 = 3\cdot 80 = 2\cdot 120$.
(Fortunately the two computations agree!)

ii) This follows from parts (ii) of the Corollaries to Theorems~1
and~2, which identify both $G$-covers with~$\SS$.~~$\Box$

The reader may have already surmised that the distinction between
even and odd theta characteristics is equivalent to the unitary
structure on $J(C')[2]$.  Indeed a theta characteristic $K_0+T$\/
is even according as $\upsilon(T)=0$ or~1.  This can be proved
as follows.  Let $q$ be the affine-quadratic form on the
theta characteristics $\Theta$ which is $0$ on even and $1$
on odd $\Theta$'s.  It is known~\cite{Mumford} that $q$
is compatible with the Weil pairing: if $T,T'$ are any 2-torsion
points then their Weil pairing $\langle T,T' \rangle$ (written
additively to take values in $\F_2$ rather than $\{\pm1\}$) is
\be
q(\Theta) + q(\Theta+T) + q(\Theta+T') + q(\Theta+T+T')
\label{4q}
\ee
for any choice of $\Theta$.  Now take $\Theta=K_0$ and
$T'=\vrho T$.  We have seen already that $q(K_0)=0$.
Moreover $q$ is $\vrho$-invariant.  Thus (\ref{4q})
reduces to
\be
q(K_0 + T) + q(K_0 + \vrho T) + q(K_0 + \vrho^2 T)
= 3 q(K_0+T) = q(K_0+T).
\label{3q}
\ee
Thus $q(K_0+T) = \langle T,\vrho T \rangle = \upsilon(T)$
as claimed.

{\bf 3. The lattice $\Erho$.}
Our proof of Theorem~3 does not entirely dispel the mystery
of the coincidence of the two $G$\/-covers of~$\Six$;
even the fact that $\PSp_4(\F_3)$ and $\SU_4(\F_2)$ are
isomorphic appears to emerge as an accident --- in fact
formulate our analysis so that the isomorphism
$\PSp_4(\F_3) \cong \SU_4(\F_2)$ arises as a by-product!
But there is a more satisfactory approach to this
isomorphism via the low-dimensional representations
of~$G$\/ and its double cover $2.G = \Sp_4(\F_3)$.
(The existence of $2.G$ cannot be so readily seen
from the $\SU_4(\F_2)$ model of~$G$\/: except for
in characteristic~2 have trivial Schur multiplier~\cite[p.xvi]{ATLAS}.)
A four-dimensional representation~$V$\/ of $2.G$\/ will also
figure in a third description of our $G$\/-cover $\SS$\/ of~$\Six$
which will naturally yield both the $J(C)[3]$ and $J(C')[2]$
pictures of~$\SS$\/ via the reduction of~$V$\/ mod~2 and~3
respectively.

It is known that $\Aut(G)$ contains~$G$\/ with index~$2$.  When
$G$\/ is viewed as $\PSp_4(\F_3)$, the outer automorphisms are
linear transformations that multiply the symplectic form by~$-1$,
while in the $\SU_4(\F_2)$ viewpoint they are conjugate-linear
transformations preserving the unitary form.  Thus $\Aut(G)$
appears as $\PGSp_4(\F_3)$ in characteristic~3 and as $\SigU_4(\F_2)$
in characteristic~$2$.

We take for $V$\/ one of the two irreducible representations
of~$2.G$\/ of dimension~$4$; the choice does not matter, because
these representations (which are each other's contragredient)
are exchanged by an outer automorphism of~$G$.  This representation
is defined over $\Q(\rho)$; the image of $2.G$ in $\GL(V)$
(actually $\SL(4)$), extended by the three-element group
$\langle\rho\rangle$ of $\bmu_3$-multiples of the identity,
is a complex reflection group --- \#32 in the list
of~\cite[p.301]{ST} --- generated by what the ATLAS calls
``triflections'': linear transformations of order~$3$\/
with a codimension-1 fixed subspace.
The nontrivial central element of~$G$\/ acts on~$V$\/
by multiplication by $-1$; thus the exterior square $W=\wedge^2 V$\/
is a 6-dimensional representation of~$G$.  It turns out that this
representation is irreducible and defined over $\Q$.  This time
$W$\/ is the unique representation of its dimension; thus the
action of~$G$\/ on~$W$\/ extends to $\Aut(G)$.  It turns out that
the image of $\Aut(G)$ in $\GL(W)$ is again a reflection group,
\#35 in~\cite{ST}, which is to say the Weyl group $W(E_6)$.
Kneser observed in~\cite{Kneser} that reducing the action
of~$\Aut(G)$ on the $E_6$ lattice mod~3 and~2
identifies $\Aut(G)$ with $\PGSp_4(\F_3)$ and $\SigU_4(\F_2)$
in these groups' guises as the orthogonal groups $\SO_5(\F_3)$ and
$\SO^-_6(\F_2)$.  Indeed $E_6/3E_6^*$ is and $E_6/2E_6$ are
orthogonal spaces of dimensions~5 and~6 over $\F_3$ and $\F_2$
respectively, both with actions of~$G$\/ respecting the quadratic
form; but the simple groups $\SO_5(\F_3)$ and $\SO^-_6(\F_2)$ are
barely large enough to accommodate copies of~$\Aut(G)$, so the reduction
maps from~$\Aut(G)$ to these two groups must be isomorphisms.  In
particular the isomorphism between $\SO_5(\F_3)$ and $\SO^-_6(\F_2)$
follows.  Kneser explains most of the sporadic isomorphisms between
classical simple groups in the same fashion.  But he does not
observe an alternative explanation along the same lines of the
isomorphism $\PSp_4(\F_3) \cong \SU_4(\F_2)$ using a $2.G$-lattice
$\Erho\subset V$\/; and it is this lattice, and its reductions
mod~$2$ and~$3$, that are central to our story.\footnote{
  Remarkably $G$\/ is involved in yet another complex reflection
  group, the 5-dimensional \#33 defined over $\sQ(\rho)$.  The
  reductions of this representation mod~$3$ and~$2$ again yield
  the isomorphism $\PSp_4(\sF_3) \cong \SU_4(\sF_2)$.  We do not
  pursue this here because this representation does not enter
  into our investigation of~$\SS$.\label{5dim}
  }

We call the lattice $\Erho$ because it is obtained from the
$E_8$ root lattice by choosing a \hbox{3-cycle} in $W(E_8)$ acting
on $E_8$ with trivial fixed space.  There is a unique conjugacy
class of such \hbox{3-cycles}; identifying one with $\rho$ gives $E_8$
the structure of a four-dimensional lattice $\Erho$ over $\Z(\rho)$.
Its group of automorphisms is the subgroup of $W(E_8)$ commuting with
$\rho$, which is isomorphic with $(2.G) \times \langle \rho \rangle$. 
We choose the isomorphism so that $\Erho \otimes \Q(\rho)$ is
our representation~$V$\/ of~$2.G$\/ rather than its contragredient.
We can also describe $\Erho$ explicitly as the sublattice of
$\Z[\rho]^4$ consisting of vectors congruent mod $(\sqrt{-3}\,)$
to a linear combination of $(1,1,1,0)$ and $(1,-1,0,1)$, i.e.\
to a vector in the ``tetracode'' in $\F_3^4$.  (That is, $\Erho$
is the $\Z[\rho]$-lattice obtained from the tetracode by
``Construction~$\rm A_c$'' \cite[p.200, Example~11b]{SPLAG}.)
The inner product $(\vec z,\vec z')$ of $\vec z = (z_1,z_2,z_3,z_4)$
with $\vec z' = (z'_1,z'_2,z'_3,z'_4)$ is
$\frac23\sum_{j=1}^4 z\0_j \bar z'_j$; if $\vec z,\vec z'\in\Erho$
then $(\vec z,\vec z') \in (2/\sqrt{-3}\,) \Z[\rho]$.
The roots (nonzero vectors of minimal norm) are the
$2\cdot3\cdot4=24$ multiples of unit vectors by $\pm\bmu_3\sqrt{-3}$,
and the $3^3 (3^2-1) = 216$ minimal lifts of the 8 nonzero vectors
of the tetracode.  (Each such vector has 1 zero and 3 nonzero
coordinates; in a minimal lift, 0 lifts to~0, and $\pm1$ lifts
to one of the three choices $\pm\bmu_3$.)  This adds up to the
familiar count of $24+216=240$ roots.  For each root~$r$\/ we obtain
a triflection $x \mapsto x + \frac12(\rho-1) (x,r) r$, and these
triflections generate $\Aut(\Erho) = (2.G) \times \langle \rho \rangle$.
Reducing mod~2, we obtain an $\F_4$-vector space $\Erho/2\Erho$ of
dimension~$4$ with a Hermitian form $\frac12(\cdot,\cdot) \bmod 2$.
Thus the group $\Aut_1(\Erho)$ of linear isometries of~$\Erho$
of determinant~$1$ maps to a subgroup of $\SU_4(\F_2)$; the kernel
consists of the identity and multiplication by~$-1$, so we find
$\Aut_1(G)/\{\pm1\}$ as a subgroup of~$\SU_4(\F_2)$.
Reducing mod~$(\sqrt{-3}\,)$, we obtain an $\F_3$-vector
space $\Erho/\sqrt{-3}\Erho$ of dimension~$4$ with a symplectic form
$\sqrt{-3}(\cdot,\cdot) \bmod (\sqrt{-3}\,)$.  Thus $\Aut_1(\Erho)$
maps to $\Sp_4(\F_3)$, this time with trivial kernel (since the
scalar multiplications by $\rho,\bar\rho$ are excluded from $\Aut_1$).
We thus find $\Aut_1(G)/\{\pm1\}$ as a subgroup of~$\PSp_4(\F_3)$.
Since $\SU_4(\F_2)$ and $\PSp_4(\F_3)$ are barely large enough to
accommodate $\Aut_1(\Erho)/\{\pm1\}$, these two inclusions must be
isomorphisms.  We have thus explained the isomorphism of
$\SU_4(\F_2)$ with $\PSp_4(\F_3)$ by regarding these two groups
as the mod-2 and mod-3 manifestation of the characteristic-zero
object $(\Erho,2.G)$.
This explanation extends to $\Aut(G)$: the outer automorphisms
of~$2.G$ act on~$\Erho$ by conjugate-linear isometries such as
$(z_1,z_2,z_3,z_4) \mapsto (\bar z_1,\bar z_2,\bar z_3,\bar z_4)$.
Such maps descend to conjugate-linear automorphisms of $\Erho/2\Erho$,
and thus extend $\SU_4(\F_2)$ to $\SigU_4(\F_2)$; on
$\Erho/\sqrt{-3}\Erho$, they act linearly but reverse the
symplectic pairing, thus extending $\Sp_4(\F_3)$ to $\GSp_4(\F_3)$
and $\PSp_4(\F_3)$ to $\PGSp_4(\F_3)$.  Therefore
$\SigU_4(\F_2) \cong \Aut(G) \cong \PGSp_4(\F_3)$.

We next consider the reduction of the $240$ roots mod~$2$
and~$\sqrt{-3}$.  Since $\frac12(r,r)=1$ for every root~$r$,
all the roots belong to one of the 120 odd classes in $\Erho/2\Erho$;
clearly $r$ and $-r$ are in the same class, and conversely if $r,r'$
are roots congruent mod $2\Erho$ then $r'=\pm r$ because at least
one of $\frac12(r\pm r')$ is a lattice vector of norm $<2$.  But
there are $240=2\cdot 120$ roots; since each odd class mod $2\Erho$
accounts for at most two of them, each of the 120 odd classes must be
represented by a pair $\pm r$ of opposite roots.  Modulo $\sqrt{-3}$,
no root is congruent to zero and every root $r$ is congruent to
$r$, $\rho r$, and $\bar\rho r$, and to no other roots $r'$ lest
\be
4 = \left| \frac{r-r'}{\sqrt{-3}} \right|^2
+ \left| \frac{\rho r-r'}{\sqrt{-3}} \right|^2
+ \left| \frac{\bar\rho r-r'}{\sqrt{-3}} \right|^2
\label{Erhotrick}
\ee
be the sum of three positive even integers.  Thus each of the
$80$ nonzero classes in $\Erho/\sqrt{-3}\Erho$ contains at most
$3$ roots, and since there are $240=3\cdot 80$ roots we again
conclude that each nonzero class represents a triple $\bmu_3 r$
of roots.  Of course these counts $240 = 3\cdot 80 = 2\cdot 120$ are
highly suggestive of the counts in the first parts of the Corollaries
to Theorem~1 and~2; we make the connection in the next section.
Note that we could also have obtained these results from the
identification of~$\Aut(\Erho)/\{\pm1\}$ with $\SU_4(\F_2)$
and $\PSp_4(\F_3)$, together with the fact that the unitary and
symplectic groups act transitively on odd and nonzero vectors
respectively.

We conclude our description of the lattice $\Erho$ by noting that
it could also have been defined directly from the representation
of $2.G$\/ on~$V\!$, without recourse to the $E_8$ root lattice:
the representation is globally irreducible, and thus has
a unique $(2.G)$-stable lattice, see \cite{Gross,Thompson}.

{\bf 4. The elliptic surface $\E$ and its \MW\ lattice.}
We next identify $\Erho$ with the group of $\bar k(t)$-rational
points on the elliptic curve
\be
\E: y^2 = x^3 + S(t)
\label{Edef}
\ee
for each sextic $S\in k[t]$ without repeated roots.  To do this
we must show not only that $\E(\bar k(t))$ is a free abelian
group of rank~$8$ but also specify an action of~$\rho$ and
a quadratic form and an action of~$\rho$ on that group.
We will have $\rho$ act by {\em complex multiplication}\/~(CM),
and take for the quadratic form the {\em canonical height}~$\Ht$\/
on the points of an elliptic curve over a function field.
We describe these extra structures (CM and $\Ht$\/) on
$\E(\bar k(t))$ in turn, and then show that the resulting
$\Z[\rho]$-lattice is isometric with~$\Erho$.

That $\E$\/ has complex multiplication by $\Z[\rho]$ means
that $\Z[\rho]$ acts on~$\E$\/ by endomorphisms (a.k.a.\ isogenies,
i.e.\ algebraic maps commuting with the group law).  Giving $\E$\/
a $\Z[\rho]$ action means exhibiting an endomorphism $\rho$
satisfying the minimal equation $\rho^2+\rho+1=0$, with addition
defined via the group law on~$\E$.  Such an endomorphism is
$\rho: (x,y)\mapsto(\rho x,y)$: the three images $(x,y)$, $(\rho x,y)$,
$(\bar\rho x,y)$ of a generic point $(x,y)$ on~$\E$\/ under
$1,\rho,\rho^2$ have the same $y$-coordinate, and thus are
the intersection of~$\E$\/ with a line and sum to zero
in the group law.  We chose $(x,y)\mapsto(\rho x,y)$ rather
than $(x,y)\mapsto(\bar\rho x,y)$ so that $\rho$ multiplies
by $\rho$ the invariant differential $\om=dx/y$ on~$\E$\/; by
linearity it follows that $\phi^*\om = \om$ for every
endomorphism $\phi\in\Z[\rho]$.  We note for future use that
any $\phi\in\Z[\rho]$, considered as a map $\phi:\E\ra\E$,
has degree $\phi\bar\phi$.  Every elliptic curve in
characteristic other than~2 or~3 that has CM by $\Z[\rho]$
can be written as $y^2=x^3+a_6$ for some nonzero $a_6$, and
so becomes isomorphic with
\be
E_0: Y^2 = X^3 + 1
\label{E0}
\ee
once we extract a sixth root of~$a_6$.  Thus in our case $\E$\/
becomes isomorphic with the constant curve~$E_0$ when we extend
$k(t)$ to the function field of the curve
\be
C'': w^6 = S(t).
\label{C''}
\ee
This is a curve of genus~10 (a smooth plane sextic)
which cyclically covers both~$C$\/ and~$C'$ in degrees~$3$
and~$2$ respectively.  The identification of~$\E$\/ with~$E_0$
over $k(C'')$ takes $(x,y)$ to $(w^{-2}x,w^{-3}y)$.
Thus the $k(t)$-rational points of~$\E$\/ are identified with
the subgroup of $\E_0(k(C''))$ consisting of the $k(C'')$-rational
points~$P$\/ of~$E_0$ whose image under the generator of
$\Gal(k(C'')/k(t))$ taking $w$ to~$-\rho w$\/ is $-\rho P$.
Note that $\E_0(k(C''))$ may be regarded as the set of maps
$P: C''\ra E_0$, which inherits an abelian group structure
from the group law of~$E_0$; in this viewpoint $\E(k(t))$
consists of those maps~$P$\/ for which the diagram
\be
\begin{array}{lcl}
C'' \!\! & \stackrel{P}{\longrightarrow} \!\! & E_0 \\
\,\downarrow & & \,\downarrow(-\rho)\\
C'' \!\! & \stackrel{P}{\longrightarrow} \!\! & E_0
\end{array}
\label{commdiag}
\ee
commutes using $(t,w)\mapsto(t,-\rho w)$ for the arrow $C''\ra C''$.
Likewise we may regard the $\bar k(t)$-rational points of~$\E$\/ as
maps from $C''$ to $E_0$ defined over~$\bar k$\/ for which that
diagram commutes.

We next review the canonical height on~$\E$\/; see
\cite[Ch.VIII\,\S9 (pp.227--233)]{Silverman} for the canonical
height on a general elliptic curve over a global field,
and \cite{CZ,Shioda,E:MW} for the present case of an
elliptic curve, especially one of constant $j$-invariant,
over a function field.

Let $E$\/ be any elliptic curve over a function field~$F=k(C_0)$,
with $C_0$ a curve over the field of constants~$k$.  Choose
a Weierstrass model for~$E$, with coordinates $x,y$.  The
{\em \naive\ height}\/ on~$E$\/ is a function~$h$\/ on its group
$E(F)$ of rational points whose value at a point measures the point's
complexity: the zero point has \naive\ height~$0$, and a nonzero point
$(x,y)$ has \naive\ height
\be
h(x,y) := \max(\deg x, \frac23 \deg y).
\label{naive}
\ee
(Here ``$\deg$'' is the degree of an element of the function field
considered as a rational map to~$\P^1$.  Note that $h$\/ depends
on the choice of Weierstrass model for~$E$, though the \naive\ heights
associated with different models differ only by $O(1)$.  The
{\em canonical} or {\em N\'eron-Tate height}\/ on the curve
does not depend on the choice of model.  It is a $\Q$-valued quadratic
form $\Ht$\/ on $E(F)$, positive-definite on $E(F)\otimes_\subZ\R$,
which is within $\pm O(1)$ of the \naive\ height~$h$.  These
two properties uniquely characterize~$\Ht$.  It follows that if
$C_1$ is a curve covering $C_0$ with degree~$d$, and
$F'/F$\/ the corresponding field extension with $F'=k(C_1)$
and $[F':F]=d$, then for any $P\in E(F)$ its height $\Ht_{F'}(P)$
as a point of $E(F')$ is $d$\/ times its height $\Ht_F(P)$ as a point
of~$E(F)$.  Indeed, from the definition~(\ref{naive}) it is clear
that the \naive\ heights $h_F,h_{F'}$ satisfy $h_{F'} = d h_F$,
so the restriction of $d^{-1}\Ht_{F'}$ to $E(F)$ satisfies both
criteria for $\Ht_F$.

There are other equivalent characterizations
more suitable for computing~$\Ht$.  Tate showed that for any
endomorphism $\phi:E\ra E$\/ of degree $>1$ the canonical height is
the unique function $\Ht: E(F)\ra\R$ such that $\Ht=h+O(1)$ and
$\Ht(\alpha P)=(\deg\alpha)\Ht(P)$ for every $P\in E(F)$.
(Usually one takes for $\phi$ the multiplication-by-$n$ map
for some small $n>1$ such as $n=2$, but we'll be able to exploit
$\phi\notin\Z$ as well.) N\'eron described $\Ht$ in terms
of intersection theory on~$E$\/ considered as an elliptic surface
over the constant field of~$F$, which yields a formula for~$\Ht(P)-h(P)$
as a finite sum of terms each depending on the reduction of~$P$\/
at one of the places where $E$\/ has bad reduction.  From either
characterization it follows that if $E$\/ is a curve defined
over~$k$, so points $P\in E(F)$ are equivalent to rational maps
$P: C_0\ra E$, then the height of such a point equals twice
the degree of~$P$\/; that is, the canonical height equals the
\naive\ height.  (Since $x,y$ are functions of degree~$2,3$
on~$E$, the definition~(\ref{naive}) gives $h(P)=2\deg(P)$.)
In Tate's viewpoint we see this by noting that $h$\/ already
satisfies the condition $h(\alpha P)=(\deg\alpha)h(P)$, because
the degree of rational maps is multiplicative under composition.
In N\'eron's approach this follows because there are no places
of bad reduction and thus no contributions to the sum for $\Ht-h$.

In our case (\ref{Edef}) of the curve $\E$\/ over~$k(t)$ or $\bar k(t)$
we likewise show that $\Ht=h$\/:

{\bf Proposition.} {\sl
The canonical height of any nonzero $P\in\E(\bar k(t))$ with coordinates
$(x,y)$ is given by $\Ht((x,y)) = \max(\deg x, \frac23 \deg y)$;
that is, the \naive\ and canonical heights on $\E(\bar k(t))$ are equal.
}

As in~\cite{E:MW} we can give several proofs: one using the
relationship with maps $C''\ra E_1$, one using Tate's characterization,
and one using N\'eron's formula.  We give the first and third proof,
leaving the second as an exercise; each of the three introduces ideas
that will later figure in our determination of $\E(\bar k(t))$
and its connection with~$\SS$.  The following observations regarding
the \naive\ height will be helpful throughout.
Let $P:(x(t),y(t))$ be a nonzero point on~$\E(\bar k(t))$.  Let
$\tau_r\in\P^1(\bar k)$ be the points at which
either $x$ or~$y$ has a pole, with $t_0=\infty$.  For $j>0$,
since $S=y^2-x^3$ is regular at $t_j\in\bar k$, the functions
$x^3,y^2$ must have poles of the same order at $\tau_r$, so $x,y$ 
ave poles of order $2n_r,3n_r$ at $\tau_r$ for some integer~$n_r$.
At $t_0=\infty$, $y^2-x^3$ has a pole of order~$5$ or~$6$,
so at least one of $x^3,y^2$ has a pole of order 6, and if one
has a pole of order $>6$ then both do and the pole orders are equal.
In this last case $x,y$ have poles of order $2n_0,3n_0$ at $t=\infty$
for some integer $n_0>1$; otherwise we define $n_0=1$.  Then
\be
h(P) = \max(\deg x, \frac23 \deg y) = 2\sum_r n_r.
\label{hsum}
\ee
We proceed to our proofs of the Proposition:

{\sl Proof 1}\/:  Extending $k(t)$ to $k(C'')$ we have identified
$\E$\/ with the constant curve $E_0$ and associated to $P:(x,y)$
the rational map $(t,w) \mapsto (w^{-2}x(t), w^{-3}y(t))$ from
$C''$ to~$E_0$.  We calculate the degree of this map by finding
the preimages of the origin of~$E_0$ and their multiplicities.
We find that the preimages are the points with $t=\tau_r$
for some~$r$, with multiplicity $6n_r$ if $\tau_r=t_j$ for some~$j$
and $n_r$ if not.  In either case $\tau_r$ contributes $6n_r$ to
the total, because there are six points of~$C''$ with $t=\tau_r$
unless $\tau_r=t_j$ in which case there is only one.  Thus the
degree of the map is $6\sum_r n_r = 3 h(P)$, and its height as a
point of $\E(k(C''))$ is twice that degree, or $6 h(P)$.
Since $[k(C''):k(t)]=6$ we divide by~$6$ to obtain the claimed
formula $h(P) = \max(\deg x, \frac23 \deg y)$ for the height of~$P$\/
as a point of $\E(k(t))$.

{\sl Proof 2}\/ (sketch): Use Tate's characterization,
taking for~$\phi$ the isogeny
\be
\sqrt{-3} = \rho - \bar\rho:
(x,y) \mapsto \bigl(
\frac{3x^3-4y^2}{3x^2}, \sqrt{-3} \, \frac{9x^3 y - 8y^3}{9x^3}
\bigr)
\label{phi3}
\ee
or
\be
2: (x,y) \mapsto \bigl(
\frac{9x^4-8xy^2}{4y^2}, -\frac{27 x^6 - 9 x^3 y^2 + y^4}{8y^3}
\bigr) .
\label{phi2}
\ee
That is, show that these multiply the \naive\ height by~$3$
and~$4$ respectively.  This can be done directly by counting
poles of the coordinates of $\sqrt{-3} \cdot (x,y)$ and
$2 \cdot (x,y)$.  (One can also invoke the computation in
Proof~1 together with the fact that $\sqrt{-3}$ and~$2$,
considered as isogenies of~$E_0$, have degrees~$3$ and~$4$
respectively.)

{\sl Proof 3}\/: In N\'eron's approach the terms in the formula for
$\Ht-h$\/ are indexed by places where $E$\/ reduces to a decomposable
curve.  We will show that $\E$\/ as no such places and thus that
$\Ht-h=0$.  Specifically, we'll show that the singular reductions
of~$\E$\/ are six cubics $y^2=x^3$, with Kodaira type~II indicating a
cusp singularity at the origin but only one component, at the roots
$t_j$ of~$S(t)$.  This is clear for finite $t$, where $\E$\/ is smooth
as a surface (since the $t_j$ are distinct)
so the reduction of the N\'eron model is obtained by just
specializing~$t$, with no further blowing up necessary.
At $t=\infty$ the coefficient $S(t)$ has a pole, so we change
coordinates to $(x',y')=(x/t^2,y/t^3)$ satisfying
${y'}^2 = {x'}^3 + (S(t)/t^6)$, which at $t=\infty$ is again smooth
as a surface and so reduces to a cuspidal cubic if some $t_j=\infty$
and to an elliptic curve if not.~~$\Box$

{\sl Corollary}: i) The canonical height of every point
$P\in\E(\bar k(t))$ is an even integer, which is positive
unless $P$\/ is the zero point.\\
ii) The points $P$\/ of height~$2$ are those coming from
minimal solutions of $y^2=x^3+S(t)$.
iii) The torsion subgroup of $\E(\bar k(t))$ is trivial.

{\sl Proof}\/: i) By the Proposition it is enough to prove
this for the \naive\ height.  The height of the zero point is~$0$.
In (\ref{hsum}) we exhibited the \naive\ height
of a nonzero point~$P$\/ as a sum of nonnegative even integers
$2n_r$, including the positive $2n_0$.  Thus $h(P)$, and so also
$\Ht(P)$, is a positive even integer as claimed.

ii) If $\sum_r n_r = 2$ then $n_0=1$ and all other $n_r=0$.
Thus $x,y$ have no poles at finite~$t$, and have poles of
orders at most $2,3$ at $t=\infty$.  Conversely it is clear
that if $x,y$ are polynomials of degrees at most $2,3$ with
$y^2=x^3+S$\/ then $P:(x,y)$ is a rational point on~$\E$\/
with $h(P)=2$, and thus by the Proposition also with $\Ht(P)=2$.

iii) If $P$\/ is a torsion point then $nP=0$ for some $n>1$, so
$\Ht(P)=n^{-2}\Ht(nP)=0$.  By part (i) it follows that $P=0$.~~$\Box$

On general principles, the group of $\bar k(t)$-rational points
of~$\E$\/ is finitely generated.  The Proposition and its Corollary
tell us that this group is torsion-free and has the structure of an
even $\Z[\rho]$-lattice $L_S$, a.k.a.\ the {\em \MW\ lattice}\/
of~$\E$, in the unitary space
$(\E(\bar k(t))) \otimes_{\sZ[\rho]} \C = L_S \otimes_{\sZ[\rho]} \C$,
whose minimal nonzero points are the minimal solutions of
$y^2=x^3+S(t)$.  We already know that there are 240 such, so
the following theorem should come as no surprise:

{\bf Theorem 4.} {\sl For every sextic polynomial $S(t)$ without
repeated roots, the lattice $L_S$ is isomorphic with $\Erho$.
}

{\sl Proof}\/:  It is enough to prove that $L_S$ is isomorphic
with $E_8$ as a $\Z$-lattice, because $\Erho$ is $E_8$ with
its unique structure as a $\Z[\rho]$-lattice.  We determine
$L_S$ as a $\Z$-lattice by relating it with the N\'eron-Severi
group $\NS(\E)$ with $\E$\/ considered as a surface over~$\bar k$.
To do this we need a model of~$\E$ smooth everywhere, including
$t=\infty$.  But we have obtained such a model in the course
of proving the $\Ht=h$\/ proposition.  This model is the union
of the two open subsets $t\neq\infty$ and $t\neq0$, each of
which is represented as a surface in $\Aff^{\!1}\times\P^2$: for
$t\neq\infty$, using coordinates $t,\;(x:y:1)$; for $t\neq0$,
using $1/t,\;(x/t^2:y/t^3:1)$.  We saw already that this surface
has no reducible fibers and that its only singular fibers are
reductions of Type~II at the six points $t=t_j$.  It follows
as in~\cite{Shioda} that $\E$ is a $\bar k$\/-rational elliptic
surface, so that its N\'eron-Severi group is the even unimodular
lattice of signature $(1,9)$, of which the zero section and
the fibers contribute a unimodular sublattice of signature
$(1,1)$ (``hyperbolic plane''), and therefore the \MW\ lattice is even
unimodular of rank~$8$.  Thus $L_S\cong E_8$ as claimed.~~$\Box$

{\sl Remark}\/: We could also have proved that $L_S$ has
$\Z[\rho]$-rank 4 using the description of $L_S$ as maps to $E_0$
from $C''$, or equivalently with its Jacobian $J(C'')$, that make
(\ref{commdiag}) commute.  Now the for $0<s<6$ the
$(-\rho)^s$-eigenspace of the action of $(t,w)\mapsto(t,-\rho w)$
on $H^1(C'')$ consists of the differentials $P(t)\,dt/w^{6-s}$
with $P\in\C[t]$ of degree at most $4-s$, and thus has dimension 
$5-s$.  For $s=3$ these are the differentials pulled back from
$H^1(C)$ via the triple cover $C''\ra C$\/; for $s=2$ and $s=4$
they are the pullbacks of $H^1(C')$ via the double cover $C''\ra C'$.
Thus $L_S$ consists of the annihilator in $\Hom(J(C''),E_0)$ of
the images of $J(C),J(C')$ in $J(C'')$.  The quotient of $J(C'')$
by the abelian variety generated by $J(C)$ and $J(C')$ is an
abelian variety $J_0(C'')$ of dimension $g(C'')-g(C)-g(C')=10-4-2=4$
equipped with an endomorphism inherited from $(t,w)\mapsto(t,-\rho w)$
that multiplies all its holomorphic differentials by $-\rho$.  Thus
$J_0(C'')$ is an abelian fourfold isogenous with $E_0^4$, and the rank
of $\Hom(J_0(C''),E_0)$ equals that of $\Hom(E_0^4,E_0\0)=\Z[\rho]^4$.
Moreover, at least in characteristic~zero we see that $J_0(C'')$,
and thus the lattice $L_S$, must be independent of~$S$, because all
components of the moduli space of fourfolds isogenous with $E_0^4$
have genus~zero.  But it seems much harder to prove that
$L_S\cong\Erho$ using this approach.

By Theorem~4 the lattice of $\bar k(t)$-points on~$\E$\/ is independent
of the choice of~$S$.  However, the field extension~$k_S$ of~$k$\/
needed to define the points of $\E(\bar k(t))$ does depend on~$S$.
It is clear that this extension is finite, because $\E(\bar k(t))$ is
generated by finitely many points.  Moreover it is a normal extension,
because Galois conjugation in $k_S$ permutes the solutions of the
equation (\ref{Surface}) defined over~$k$.  Furthermore, $\Gal(k_S/k)$
respects the group law and canonical height on~$L_S$.  Thus we may
regard $L_S$ as a $\Gal(k_S/k)$ module, and obtain a map from
$\Gal(k_S/k)$ (or even from $\Gal(\bar k/k)$) to
$\Aut(L_S)\cong\Aut(\Erho)=\bmu_3\times(2.G)$.
We have seen already that, at least modulo the center $\bmu_6$
of $\Aut(\Erho)$, this Galois representation depends
only on the $\PGL_2(k)$-orbit of~$S$\/; in the $L_S$ viewpoint,
this is because if $S,S'$\/ are sextics equivalent under some 
$g\in\PGL_2(k)$ then $g:\P^1\ra\P^1$ lifts to an isomorphism from
$\E_S$ to $\E_{S'}$, using the models of these elliptic surfaces
constructed during the proof of Theorem~4.
Our work in \S2 shows that the representation
$\Gal(k_S/k)\ra G$\/ is equivalent to the Galois representations
on $J(C)[3]/\{\pm1\}$ and $J(C')[2]$, and is surjective onto~$G$\/
for generic~$S$.  But the curves $C,C'$ also figure in the arithmetic
of~$\E$\/ as an elliptic curve over~$k(t)$: their function fields are
the extensions of $k(t)$ generated by the $(\!\sqrt{-3})$-
and $2$-torsion points of the curve.  This is clear from
our formulas (\ref{phi3},\ref{phi2}) for the isogenies
$\sqrt{-3},2:\E\ra\E$\/: the nontrivial $(\sqrt{-3}\,)$-torsion
points are those with $x=0$ and thus $y=\pm\sqrt S = \pm u$,
while the nontrivial 2-torsion points have $y=0$ and
$x=-\root3\of S = -\rho^s v$.  The 3- and 2-torsion groups
of $J(C)$ and $J(C')$ then arise naturally in the {\em descent}\/
on $\E$\/ via the isogenies $\sqrt{-3}$ and~$2$, and allow us
to extend the maps of Theorems~1 and~2 from the 240 roots to
$J(C)[3]$ and $J'(C)[2]$ to isomorphisms
$L_S/\sqrt{-3}\, L_S \stackrel\sim\ra J(C)[3]$,
$L_S/2L_S \stackrel\sim\ra J(C')[2]$.  We thus obtain
our final explanation of the equivalence between the
$G$-covers of $\Six$ coming from $J(C)[3]$ and $J'(C)[2]$:
as with the mod-3 and mod-2 manifestations of~$G$\/ itself,
the two constructions of the $G$\/-cover $\SS\bigl/\Six\bigr.$ 
are now revealed as the mod-$(\sqrt{-3}\,)$ and mod-2 manifestations
of the projective representation of $\Gal(\bar k/k)$ on $L_S$.
We next exhibit these descent maps; see \cite[Ch.VIII]{Silverman}
for the general theory of descent of which they are two special cases.

We deal first with the 3-isogeny $\sqrt{-3}$.  We observed that its
kernel is generated by the 3-torsion point $(0,u)$ with $u^2=S(t)$.
This the function $y-u\in(k(C))(\E)$ is a {\em Weil function}\/
on~$\E$\/: its divisor is $3((0,u)) - 3({\bf 0})$.  Moreover $y-u$
is locally a cube at its triple pole.  Evaluation of $y-u$ mod cubes
at points other than ${\bf 0}$ and $(0,u)$ then extends to a
homomorphism from $\E(k(C))$ to $k(C)^*/{k(C)^*}^3$ whose kernel
is $\sqrt{-3}\, \E(k(C))$.  We claim that
$\sqrt{-3}\, \E(k(C)) \cap \E(k(t)) = \sqrt{-3}\, \E(k(t))$:
otherwise we have $Q=\sqrt{-3} Q_1$ with $Q\in\E(k(t))$ but
$Q_1$ defined not over $k(t)$ but only over $k(C)$;
but then the other preimages of~$Q$\/ under $\sqrt{-3}$ would be
the translates of $Q_1$ by $(0,\pm u)$, which are also defined
over~$k(C)$, and we would obtain a cubic extension of~$k(t)$ split
by the quadratic extension $k(C)/k(t)$, which is impossible.  Thus
the restriction of our homomorphism from $\E(k(C))$ to $L_S=\E(k(t))$
yields an injection from $L_S / \sqrt{-3}\, L(S)$ to
$k(C)^*/{k(C)^*}^3$.  Since the divisor of $y-u$ on~$\E$\/ is
divisible by~$3$, and $\E$ has no reducible fibers, the divisor of
$y(t)-u$ on~$C$\/ is also divisible by~$3$; call it $3D$.  (This
generalizes the argument of Theorem~1, and can also be shown as was
done there without explicitly invoking the fibers of~$\E$, though the
key point that $S$\/ has distinct roots amounts to the same thing;
in the context of descent the condition that the zero or pole
multiplicity of $y(t)-u$ be a multiple of 3 at each point $(t_0,u_0)$
of~$C$\/ is the condition of local triviality at $t_0$ of a
principal homogeneous space for the isogeny $\sqrt{-3}$.)
Then $[D]$ is a 3-torsion point on $J(C)$, and the map
$(x,y)\mapsto[D]$ is an injection from $L_S / \sqrt{-3}\, L(S)$
to $J(C)[3]$.  It is clear that the restriction of this descent map
to a root agrees with the nonzero element of $J(C)[3]$ associated
with the root in Theorem~1.

The kernel of the isogeny~2 is not cyclic, though its nonzero
elements $(0,-\root 3\of S)=(0,-v)$, $(0,-\rho v)$, $(0,-\bar\rho v)$
are permuted by~$\rho$.  The Weil functions for these points
are $x+v$, $x+\rho v$, $x+\bar\rho v$; their product
$x^3+v^3=x^3+S(t)=y^2$ is a square.  Thus evaluation of any
two of these functions, say $x+v$ and $x+\rho v$, at points
outside $\E[2]$ extends to a homomorphism from $\E(k(C'))$ to
$\bigl(k(C')^*/{k(C')^*}^2\bigr)^2$ whose kernel is $2 \E(k(C'))$.
Better yet, we can evaluate all three functions and land
in $\bigl(k(C')^*/{k(C')^*}^2\bigr)^3_0$, by which we mean the
subgroup of $\bigl(k(C')^*/{k(C')^*}^2\bigr)^3$ consisting of
triples whose product is the identity.
We again claim that $2 \E(k(C')) \cap \E(k(t)) = 2 \E(k(t))$:
otherwise we have $Q=2 Q_1$ with $Q\in\E(k(t))$ but
$Q_1$ defined not over $k(t)$ but only over $k(C')$;
but then the other preimages of~$Q$\/ under multiplication by~$2$
would be the translates of $Q_1$ by $(0,-\rho^s v)$, which are also
defined over~$k(C')$, and we would obtain a biquadratic extension
of~$k(t)$ split by the cubic extension $k(C')/k(t)$, which is
impossible.  Thus the restriction of our homomorphism from
$\E(k(C'))$ to $L_S=\E(k(t))$ yields an injection from $L_S / 2L_S$
to $\bigl(k(C')^*/{k(C')^*}^2\bigr)^3_0$.  Moreover this action is
consistent with the action of~$\rho$: substituting $\rho x$ for~$x$
cyclically permutes $x+v$, $x+\rho v$, $x+\bar\rho v$\/ modulo squares
(indeed modulo $\bmu_3$), and thus agrees with cyclic permutation of
the three coordinates of $\bigl(k(C')^*/{k(C')^*}^2\bigr)^3_0$.
Thus the image of our map is contained in the subgroup of
$\bigl(k(C')^*/{k(C')^*}^2\bigr)^3_0$ consisting of triples
of the form $(f,\bar\rho^* f,\rho^* f)$, and we may consider
it as an injection of $\F_4$-modules
$L_S/2L_S \hookrightarrow (k(C')^*/{k(C')^*}^2)$.
As in the previous paragraph we see that each of the three
functions $x+\rho^s v$ has divisor divisible by~$2$, and thus identify
the image of our map in $(k(C')^*/{k(C')^*}^2)$ with $J(C')[2]$.
Once more this is the descent map on~$\E$\/ over~$k(t)$, this time
for the multiplication-by-2 isogeny, and it is clear that the
restriction of this descent map to a root recovers the odd
element of $J(C')[2]$ associated with the root in Theorem~2.

To summarize: we have shown

{\bf Theorem 5.} {\sl
The maps from the 240 representations of~$S(t)$ as $y^2-x^3$
to $J(C)[3]$ and $J(C')[2]$ obtained in Theorems~1 and~2 extend
to surjective homomorphisms from $L_S$ to $J(C)[3]$ and $J(C')[2]$
whose kernels are $\sqrt{-3}\, L_S$ and $2L_S$ respectively;
these homomorphisms are the descent maps for the isogenies
$\sqrt{-3}$ and~$2$ on~$\E$.  The moduli space $\SS$ is the
moduli space for sextics $S(t)$ equipped with an isomorphism
of $L_S$ with $\Erho$; reducing $\Erho$ mod~\/$\sqrt{-3}$ and\/~$2$
recovers the $J(C)[3]$ and $J(C')[2]$ constructions of~$\SS$.
}

{\bf 5. Complements and coming attractions.}

{\em Orbits of the point stabilizer and the height pairing.}
The simple group~$G$\/ acts transitively on the 40 representations
of~$S$\/ as the difference between the cube of a quadratic
and the square of a cubic polynomial; but its action is not
doubly transitive: the orbits of a point stabilizer have sizes
1, 12, 27.  It is easy to see this from either the mod-3 or
the mod-2 definitions of~$G$\/: in each case the choice of
a line in a symplectic $\F_{\!3}^4$ or an odd line in a unitary
$\F_{\!4}^4$ divides the remaining 39 such lines into 12 orthogonal
to the first and 27 not orthogonal to it.  But Clebsch already
obtained the $12+27$ partition in~\cite{Clebsch}, though he
knew neither the Weil pairing on $J(C)[3]$ nor the relevance
of~$C'$.  How, then, could he distinguish the orthogonal from
the non-orthogonal pairs?  In effect he did it using the height
pairing!  Given a root $r\in\Erho$, the $6\cdot 39$ roots not
proportional to~$r$ consist of $6\cdot 12$ orthogonal to~$r$
and $6\cdot 27$ not orthogonal to it; if $r,r'$ are roots neither
proportional nor orthogonal to each other then their inner product
$(r,r')$ is one of $\bmu_6 \cdot 2/\sqrt{-3}$.  Now $(r,r')$ is
determined by the norms of $r+r'$, $r+\rho r'$, and $r+\rho^2 r'$,
which are even integers whose sum is~12 (cf.~(\ref{Erhotrick})).
If $r,r'$ are not proportional then these norms are positive and $<8$;
clearly $(r,r')=0$ if and only if each norm equals~$4$, and so
$(r,r')\neq 0$ if and only if the norms are $2,4,6$ in some order.
In the $\Erho$ picture, $r,r'$ are minimal solutions $(x(t),y(t))$
and $(x'(t),y'(t))$ of $y^2=x^3+S$, and $r+\rho^s r'$ is a solution
of $y^2=x^3+S$\/ whose coordinates are obtained from $(x,y)$ and
$(\rho^s x', y')$ by the group law on~$\E$.  The degrees of these
coordinates as functions of~$t$\/ then determine the norm of
$r+\rho^s r'$ via our ``$\Ht=h$\/'' Proposition.  The coordinates
of $r+\rho^s r'$ involve the slope quotient $(y-y')/(x-\rho^s x')$,
and thus hinge on the distribution of the zeros of $y-y'$ among
the zeros of $x-\rho^s x'$.  (Since $y^2-{y'}^2 = x^3-{x'}^3$,
the zeros of the three quadratic polynomials $x-\rho^s x'$ are 
the same as the zeros of the two cubics $y\pm y'$.)  We find that
$(r,r')=0$ exactly when each of $x-\rho^s x'$ contains just one of
the linear factors of $y-y'$, which is exactly Clebsch's condition
for distinguishing $6\cdot 12$ of the remaining $6\cdot 39$ minimal
solutions.

{\em Cubic curves tangent to six given concurrent lines.}
The problem of enumerating solutions of $y^2=x^3+S(t)$ has also
appeared in \cite{Tyrrell}, where they arose in yet another
geometric guise: the forty cubic curves tangent to the six lines
$t=t_j$ in the $(s,t)$ plane, concurrent at the point
$(1:s:t)=(0:1:0)$ at infinity.  The connection of these
cubics with $J(C)[3]$ is then seen as a special case of
the construction of the last section of~\cite{EEHS}.
We explain these connections next.

On a generic plane cubic~$\sfC$ in the $(s,t)$-plane, $t$\/ is a
rational function of degree~$3$, and thus has six ramified points,
each with one simple and one double preimage on~$\sfC$,
corresponding to the six lines through $(0:1:0)$ tangent to~$\sfC$.
Now fix six distinct points $t_j\in\P^1$ and consider curves~$\sfC$
of genus~1 with a degree-3 map $t:\sfC\ra\P^1$ ramified at the $t_j$.
By the Riemann existence theorem for branched covers of~$\P^1$,
such curves are in 1:1 correspondence with 6-tuples
of involutions in the symmetric group $S_3$ whose product is
the identity and which generate a transitive subgroup
of~$S_3$, modulo conjugation in the subgroup they generate.
Ignoring the transitivity condition and the conjugations we find
$3^5=243$ such 6-tuples: the first five involutions may be chosen
arbitrarily, and uniquely determine the sixth.  Of those $243$,
the transitivity condition excludes only those for which all six
involutions are the same, leaving an adjusted total of $243-3=240$.
Each of these must generate all of~$S_3$, so our final count of
triple covers of~$\P^1$ ramified at $t=t_j$ is $240/3!=40$.

We can directly relate these 40 covers with with two of our
three pictures of~$\SS$.  First, they naturally biject with
solutions of $y^2=x^3+S(t)$ in polynomials $x(t),y(t)$ of degrees
at most $2,3$, up to equivalence $(x,y)\sim(\bmu_3 x,\pm y)$.
Second, they naturally biject with pairs $\pm D$\/ of nontrivial
3-torsion points on the Jacobian of $C: u^2 = S(t)$.

The latter bijection was given in~\cite{EEHS};
more generally, for any distinct points $t_1,\ldots,t_{2g+2}$
in~$\P^1$, a bijection was constructed between genus-$(g-1)$
triple covers of~$\P^1$ simply ramified at the $t_j$ and pairs of
nontrivial 3-torsion points on the Jacobian of the genus-$g$
hyperelliptic curve $u^2=\prod_{j=1}^{2g+2} (t-t_j)$.
This curve is the discriminant of the triple cover;
its compositum (a.k.a.\ fiber product over~$\P^1$) with
the triple cover is an $S_3$ cover of the line which is
an unramified cyclic cubic cover of the hyperelliptic curve.
By geometric Kummer theory such covers of any curve correspond
bijectively with pairs of nontrivial 3-torsion points on the curve's
Jacobian.  For a hyperelliptic curve one readily recovers from its
cyclic cubic cover a triple cover of~$\P^1$ ramified at the
curve's Weierstrass points.  Given the $t_j$ there are thus
$(3^g-1)/2$ such triple covers, whose coefficients generate a field
extension of generic Galois group contained in $\PSp_{2g}(\F_3)$.
The count of $(3^g-1)/2$ can also be obtained as in the previous
paragraph by solving equations in~$S_3$.  Before~\cite{EEHS} it was
already known~\cite{Cohen} that the Galois group is $\PSp_{2g}(\F_3)$,
but only by combining a transcendental monodromy computation with
a difficult group-theoretical characterization of the permutation
representation of $\PSp_{2g}(\F_3)$ on $\P^{2g-1}(\F_3)$.
In \cite{EEHS} the bijection with 3-torsion in the Jacobian
was used to explain the symplectic structure in terms of the
Weil pairing on the 3-torsion.

On the other hand, it is noted in~\cite{Cohen} that from $\sfC$
and $t:\sfC\ra\P^1$ one may construct a solution of $y^2=x^3+S(t)$.
Any function of degree~$3$ on a curve~$\sfC$ of genus~$1$
is the quotient of two sections of a divisor~$D$\/ of degree~$3$;
by Riemann-Roch, $D$\/ has a three-dimensional space of sections,
which identifies $\sfC$ with a cubic curve in the $(s,t)$ plane
up to affine-linear transformations preserving~$t$, i.e.\ up to
transformations of the form $(t,s)\mapsto(t,as+bt+c)$.
Using these transformations we put our cubic in the form
$s^3-3x(t)s+2y(t)=0$ where the polynomials $x,y$ are of degree
at most~$2$ and~$3$ respectively.  This form is unique up to scalings
$(s,x,y)\mapsto(\lambda s, \lambda^2 x, \lambda^3 y)$.
The ramified points are the values of~$t$\/ at which the
discriminant of this cubic in~$s$\/ vanishes, i.e.\ the
roots of the sextic $y^2-x^3$.  Scaling by $\lambda$
multiplies $y^2-x^3$ by~$\lambda^6$.  Thus there is a
$\lambda$, unique up to multiplication by $\bmu_6$, that
yields $y^2-x^3=S(t)$.  Conversely, from a solution to
$y^2-x^3=S(t)$ we recover a plane cubic $\sfC:s^3+3x(t)s+2y(t)=0$
with a degree-3 cover $t:\sfC\ra\P^1$ ramified at the $t_j$, with
two solutions $(x,y)$ producing the same curve if and only if
they are equivalent under $(x,y)\sim(\bmu_3 x,\pm y)$.  At this point
in~\cite{Cohen}, Cohen cites Clebsch for the enumeration of
solutions of $y^2-x^3=S(t)$; indeed \cite{Cohen} was the paper
that alerted us to Clebsch's work on this problem.  It is easy
to see that the pair of 3-torsion points of $J(C)$ constructed from
$y^2-x^3=S(t)$ in Theorem~1 is the same as the pair obtained from
the cubic $s^3+3x(t)s+2y(t)=0$ using the construction of~\cite{EEHS}.

{\em Ground fields without roots of unity.}
We have assumed throughout that our ground
field~$k$\/ contains the cube roots of unity.  If we work over
a field, such as~$\Q$, that does not contain them but still has
characteristic other than~$2$ or~$3$, then we can still define
the moduli spaces $\Six$ and $\SS$, but the latter space's function
field must contain $\bmu_3$.  This can be seen from each of our
three pictures of~$\SS$\/: in the $J(C')[2]$ and \MW\ pictures, 
via the action of~$\rho$ on $C'$ and~$\E$\/; in the $J(C)[3]$ picture,
via the Weil pairing.  [For the latter, it might seem that the sign
ambiguity in $J(C)[3]$ may frustrate the extraction of a cube root of
unity from the level-3 structure.  But the only ambiguity is the
possibility of multiplying all of $J(C)[3]$ by~$-1$, and if $P,Q$\/
are 3-torsion points whose Weil pairing $\langle P,Q \rangle$ is~$\rho$
then $\langle -P,-Q \rangle = \rho$ also.]  The field extension
$k\bigl(\SS\bigr)$ of $k\bigl(\Six\bigr)$ then has Galois group
$\Aut(G)$, with the outer automorphism of~$G$\/ inducing the
Galois involution of $k(\rho)/k$.

{\em Twists, and $G$\/ vs.\ $\Aut(\Erho)$.}  From a sextic
$S(t)\in k[t]$ without repeated roots we find, from the fiber
of~$\SS$\/ above the associated point of~$\Six$, a normal extension
$k'$ of~$k$\/ of degree at most $25920$ and a map $\Gal(k'/k)\ra G$.
But this $k'$ is not the field of definition~$k_S$ of~$L_S$: for
each minimal solution $(x,y)$ of $y^2-x^3=S(t)$ the coordinates
of $x^3,y^2$ are contained in~$k'$, but not necessarily the
coordinates of~$x,y$ themselves.  Indeed it might seem that to obtain
$k'$ from~$k$\/ we might need to extract a different sixth root
for each of the 40 minimal representations of $S(t)$ as the difference
between a square and a cube.  But in the $\Erho$ picture it is
clear that a single cyclic sextic extension suffices to obtain
$k_S$ from~$k'$: the group $G$\/ is the quotient of $\Aut(\Erho)$
by its center~$\bmu_6$.  Note, however, that the extension $k_S/k'$
is not visible at the level of moduli spaces: multiplying $S$\/
by $\lambda\neq0$ does not change the associated point of~$\Six$
and thus gives rise to the same extension $k'$, but multiplies
each of $x^3,y^2$ by $\lambda$ and thus twists the extension $k_S/k'$
by $\root6\of\lambda$.  Clearly $k_S/k'$ is the compositum of
two extensions of $k'$, the first obtained by adjoining the
coefficients of the $x$-coordinates of the roots, the second obtained
by adjoining the coefficients of the $y$-coordinates.  The former
extension is either trivial or cyclic cubic, and the latter is
either trivial or quadratic.  Since the double cover~$2.G$ of~$G$\/
does not split, the quadratic cover, though normal over~$k$,
cannot be obtained from $k_S$ by adjoining the square root
of an element of~$k$.  On the other hand the index-2 subgroup
of $\Aut(\Erho)$ does split as $\bmu_3\times G$, so the subfield
of $k_S$ obtained by adjoining each root's $x$~coefficients is
the compositum of~$k'$ with a cyclic cubic extension of~$k$.
It might be surmised that this extension is
$k(\root3\of{c\Delta(S)})$ for some $c\in\Q^*$,
where $\Delta(S)$ is the discriminant of~$S$\/; since
$\Delta(\lambda S)=\lambda^{10}\Delta(S)$, this surmise
would behave correctly under scaling of~$S$.  In a future
paper we shall show that in fact $k(\root3\of{\Delta(S)})$
(with $c=1$) is the correct extension.

{\em Geometry of~$\SS$ and its fiber products with $\Six$ and $\Five$.}
We have given three descriptions of the moduli space $\SS$\/ but
have not described its geometry.
Burkhardt~\cite{Burkhardt} had already in effect identified $\SS$ with
an open set in his quartic hypersurface in $\P^4$, which is the zero
locus of the degree-4 invariant of the five-dimensional representation
of~$G$\/ mentioned in footnote~\ref{5dim}; see also \cite[p.190]{Hunt}.
This hypersurface is rational \cite[p.184--5]{Hunt}, a fact Hunt
attributed to Todd (1936), with an explicit birational map to $\P^3$
first given by Baker six years later.  The fiber product over~$\Six$
of~$\SS$ with $\SIX$ --- that is, the moduli space
of genus-2 curves~$C$\/ with full level-6 structure, or equivalently of
curves~$C'$ with full level-$(2\sqrt{-3})$ structure --- is an open set
in an algebraic threefold of general type~\cite{Janus}.  Recall that
this moduli space, call it $\SS'$, arose in the course of Hunt and
Weintraub's identification of the $G$\/-covers of~$\Six$ coming from
$J(C)[3]$ and $J(C')[2]$; it is an $S_6$~cover of~$\SS$.
We have investigated an intermediate, non-Galois cover of~$\Six$,
which is the fiber product of~$\SS$\/ with the moduli space of
configurations of six points of~$\P^1$ one of which is distinguished.
We may put this point at infinity, and consider the remaining five
as points on the affine line $\Aff^1$; in keeping with our earlier
notation we thus call this moduli space $\Five$.
Using the $\E$\/ picture of~$\SS$, and specializing the
formulas of~\cite{Shioda}, we have shown that the resulting 
$G$\/-cover of~$\Five$ is still rational, and is geometrically
even nicer than~$\SS$\/: it is the complement of hyperplanes
in~$\P^3=P(V)$!  The representation of $2.G$\/ on~$V$\/ yields
the action of~$G$\/ on this projective space, and the excluded
hyperplanes are the 40 orthogonal complements of the roots of~$\Erho$.
We shall show this in a future paper; the computations there will
also verify the claim in the previous paragraph that
$k(\root3\of{\Delta(S)}) \subseteq k_S$.

{\em Degenerations.}  In the same paper we shall also describe the
behavior of the configuration of the $t_j$, and thus also of
$J(C)[3]$, $J(C')[2]$ and~$\E$, as we approach the excluded
hyperplanes.  For instance, we find that as two of the five
finite $t_j$ approach one another, the point on~$\P^3$ approaches
one of the hyperplanes, while as some $t_j\ra\infty$ the point
approaches one of the 40 (projectivizations of) roots, which must
be blown up to detect the configuration of the remaining four~$t_j$.
Perhaps most strikingly, we obtain the moduli space $\A_2(3)$ of
principally polarized abelian surfaces with full level-3 structure
by blowing up $\P^3$ at the 90 lines such as $z_1=z_2-z_3=0$
that contain four of the 40 projectivized roots: $\A_2(3)$ is
the complement in that blown-up $\P^3$ of the proper transforms
of the 40 hyperplanes.

\end{document}